\documentstyle[12pt,amssymb,amstex]{amsart}

\newcommand{\bea}{\begin{eqnarray}}
\newcommand{\eea}{\end{eqnarray}}
\newcommand{\ba}{\begin{array}}
\newcommand{\ea}{\end{array}}
\newcommand{\edc}{\end{document}}
\newcommand{\bc}{\begin{center}}
\newcommand{\ec}{\end{center}}
\newcommand{\be}{\begin{equation}}
\newcommand{\ee}{\end{equation}}

\newcommand{\dsf}{\displaystyle\frac}
\def\s{\sigma}
\def\e{{\bf 1}\!\!{\rm I}}
\def\l{\lambda}
\def\z{\eta}
\def\t{\tau}

\def\a{\alpha}
\def\w{\omega}
\def\O{\Omega}

\def\K{\mathcal G}
\def\R{\mathbb{R}}
\def\b{\beta}
\def\m{\mu}

\def\L{\Lambda}
\def\G{\Gamma}
\def\C{\mathbb{C}}

\def\Z{\mathbb{Z}}

\newtheorem{theorem}{Theorem}[section]
\newtheorem{lemma}[theorem]{Lemma}

\begin{document}

\title[]{on factors associated with quantum Markov states corresponding
to nearest neighbor models on a Cayley tree}
\author{Francesco Fidaleo}
\address{Dipartimento di Matematica\\
Universit\`{a} di Roma ``Tor Vergata''\\
Via della Ricerca Scientifica, 00133 Roma, Italy,}
\email{{\tt fidaleo@@axp.mat.uniroma2.it}}
\author{Farruh Mukhamedov}
\address{Department of Mechanics and Mathematics\\
National University of Uzbekistan\\
Vuzgorodok, 700095, Tashkent, Uzbekistan,}
\email{{\tt far75m@@yandex.ru}}

\begin{abstract}
In this paper we consider nearest neighbour models where the spin
takes values in the set $\Phi=\{\z_1,\z_2,...,\z_q\}$ and is
assigned to the vertices of the Cayley tree ${\G}^k$. The
Hamiltonian is defined by some given $\l$-function. We find a
condition for the function $\l$ to determine the type of the von
Neumann algebra  generated by the GNS - construction associated
with the quantum Markov state corresponding to the unordered phase
of the $\l$-model. Also we
give some physical applications of the obtained result.
\end{abstract}

\maketitle

\section{Introduction}

It is known that in the quantum statistical mechanics concrete
systems are identified with states on corresponding algebras. In
many cases the algebra can be chosen to be a quasi-local algebra
of observables.  The states on this algebra satisfying the KMS
condition, describe equilibrium states of the quantum system.
Basically, limiting Gibbs measures of classical systems with the
finite radius of interactions are  Markov random fields (see e.g.
\cite{D},\cite{P}). In connection with this, there arises a problem to
construct analogues of non-commutative Markov chains. In \cite{A}
Accardi explored this problem, he introduced and studied
noncommutative Markov states on the algebra of quasi-local
observables which agreed with the classical Markov chains. In
\cite{AF},\cite{GZ},\cite{AFi}, modular properties of the 
non-commutative
Markov states were studied. In \cite{FNW}  Fannes, Nachtergale and
Werner showed that ground states of the valence- bond- solid modles on 
a Cayley
tree  were quantum Markov chains on the quasi-local algebra. In
the present paper we will consider Markov states associated with
nearest neighbour models on a Cayley tree. Note the investigation
of the type of quasi-free factors (i.e. factors generated by
quasi-free representations) has been an interesting problem since
the appearance of the pioneering work of Araki and Wyss \cite{AW}.
In \cite{Po} a family of representations of uniformly hyperfinite
algebras was constructed, which can be treated as a free quantum
lattice system. In that case factors corresponding to those
representations are of type III${}_{\l}$, $\l\in (0,1)$. More
general constructions  of product states were considered in
\cite{AWo}.

Observe  that the product states can be viewed as the Gibbs states
of Hamiltonian system in which interactions between particles of
the system are absent, i.e. the system is a free lattice quantum
spin system. So it is interesting to consider quantum lattice
systems with non-trivial interactions, which leads us, as
mentioned above, to consider the Markov states. Simple examples of
such systems are  the Ising and Potts models, which have been
studied in many papers (see, for example, \cite{SS}). We note that
all Gibbs states corresponding to these models are Markov random
fields. The full analysis of the type of von Neumann algebras
associated with the quantum Markov states is still an open
problem. Some particular cases  of the Markov states were
considered in \cite{FM},\cite{GN},\cite{M},\cite{MR1},\cite{MR2}.

The present paper is devoted to the type analysis of some class of
diagonal quantum Markov states, which correspond to a $\l$-model
on the Cayley tree, in which spin variables take their values in a
set $\Phi=\{\z_1,...,\z_q\}$, where $\z_k\in{\R}^{q-1}$,
$k=1,\dots,q$. Observe that the considered model generalizes a
notion of $\l$-model introduced in \cite{R}, where the spin
variables take their values $\pm 1$.

\section{Definitions and preliminary results}

The Cayley tree  $\Gamma^k$ of order $ k\geq 1 $ is an infinite
tree, i.e., a graph without cycles, such that each vertex of which
lies on $ k+1 $ edges. Let $\Gamma^k=(V, \Lambda),$  where $V$ is
the set of vertices of $ \Gamma^k$, $\Lambda$ is the set of edges
of $ \Gamma^k$. The vertices $x$ and $y$ are called {\it nearest
neighbor}, which is denoted by $l=<x,y>$ if there exists an edge
connecting them. A collection of the pairs
$<x,x_1>,...,<x_{d-1},y>$ is called {\it path} from the point $x$
to the point $y$. The distance $d(x,y), x,y\in V$, on the Cayley
tree, is the length of the shortest path from $x$ to $y$.

We set
$$ W_n=\{x\in V| d(x,x^0)=n\}, $$
$$ V_n=\cup_{m=1}^n W_m=\{x\in V| d(x,x^0)\leq n\}, $$
$$ L_n=\{l=<x,y>\in L | x,y\in V_n\}, $$
for an arbitrary point $ x^0 \in V $. Denote $|x|=d(x,x_0)$, $x\in
V$.

Denote
$$
S(x)=\{y\in W_{n+1} :  d(x,y)=1 \}, \ \ x\in W_n.
$$

This set is made of the {\it direct successors} of $x$. Observe
that any vertex $x\neq x^0$ has $k$ direct successors and $x^0$
has $k+1$.

\begin{theorem} \label{21}
There exists a one-to-one correspondence between the set $V$ of
vertices of the Cayley tree of order $k\geq 1$ and the group
$G_{k+1}$ of the free products of $k+1$ cyclic  groups of the
second order with generators $a_1,a_2,...,a_{k+1}$.
\end{theorem}

Consider a left (resp.  right) transformation shift on $G_{k+1}$
defined as follows. For $ g_0\in G_{k+1}$  we  put
$$
T_{g_0}h=g_0h \ \ (resp. \ T_{g_0}h=hg_0,), \ \  \ h\in G_{k+1}.
$$

It is easy  to see that  the set  of all left  (resp. right)
shifts on $G_{k+1}$ is  isomorphic to the group $G_{k+1}$.

Let $\Phi=\{\z_1,\z_2,...,\z_q\}$, where $\z_1,\z_2,...,\z_q$ are
vectors in ${\R}^{q-1},$ such that
\begin{equation}
\z_i\z_j = \left \{ \ba{ll}
1, \  \ \textrm{if} \ \ i=j \\
-\frac{1}{q-1}, \  \ \textrm{if} \ \  i \neq j \ea \right.\,.
\label{2.1}
\end{equation}

We consider models where the spin takes values in the set
$\Phi=\{\z_1,\z_2,...,\z_q\}$ and is assigned to the vertices of
the tree. A configuration $\s$ on $V$ is then defined as a
function $x\in V\to\s(x)\in\Phi$; the set of all configurations
coincides with $\Omega=\Phi^{\Gamma^k}$. The Hamiltonian is of an
$\l$-model form :
\begin{equation} H_{\l}(\s)=\sum_{<x,y>}\l(\s(x),\s(y);J)\,,
\label{2.2} \end{equation} where $J\in {\R}^n$ is a coupling
constant and the sum is taken over all pairs of neighboring
vertices $<x,y>$, $\s\in \Omega$. Here and below
$\l:\Phi\times\Phi\times {\R}^n\to {\R}$ is some given function.

We note that $\l$-model of this type can be considered as a
generalization of the Ising model. The Ising model corresponds to
the case $q=2$ and $\l(x,y;J)=-Jxy$.

We consider a standard $\s$-algebra ${\mathcal F}$ of subsets of
$\Omega$ generated by cylinder subsets, all probability measures
are considered on $(\Omega,{\mathcal F})$. A probability measure $\m$
is called a {\it Gibbs measure} (with Hamiltonian $H_{\l}$) if it
satisfies the DLR equation: for $n=1,2,...$ and $\s_n\in
\Phi^{V_n}$:
\begin{equation*} \m\bigg(\{\s\in \O :
\s|_{V_n}=\s_n\}\bigg)=\int_{\O}\m(d\w)\nu^{V_n}_{\w
|_{W_{n+1}}}(\s_n)\,
\end{equation*}
where $\nu^{V_n}_{\w|_{W_{n+1}}}$ is the conditional probability
$$ \nu^{V_n}_{\w |_{W_{n+1}}}(\s_n)=Z^{-1}(\w |_{W_{n+1}})\exp(-\b
H(\s_n||\w|_{W_{n+1}})). $$ where $\b>0$. Here $s_n|_{V_n}$ and
$\w|_{W_{n+1}}$ denote the restriction of $\s,\w\in \O$ to $V_n$
and $W_{n+1}$ respectively. Next, $\s_n:x\in V_n\to\s_n(x)$ is a
configuration in $V_n$ and $ H(\s_n||\w|_{W_{n+1}})$ is defined as
the sum $H(\s_n)+U(\s_n,\w|_{W_{n+1}})$ where
\begin{align*}
H(\s_n)=&\sum_{<x,y>\in L_n}\l(\s_n(x),\s_n(y);J), \\
U(\s_n,\w|_{W_{n+1}})=&\sum_{<x,y>: x\in V_n, y\in
W_{n+1}}\l(\s_n(x),\w(y);J)\,.
\end{align*}

Finally, $Z(\w |_{W_{n+1}})$ stands for the partition function in
$V_n$ with the boundary condition $\w|_{W_{n+1}}$:
\begin{equation*} Z(\w
|_{W_{n+1}})=\sum_{\tilde\s_n\in\Phi^{V_n}}\exp(-\b
H(\tilde\s_n||\w|_{W_{n+1}})\,. \end{equation*}

Since  we consider nearest neighbour interactions, the Gibbs
measures of the $\l$-model possess a Markov property: given a
configuration $\w_n$ on $W_n$, random configurations in $V_{n-1}$
and in $V\setminus V_{n+1}$ are conditionally independent. It is
known (see \cite{S}) that for any sequence $\w^{(n)}\in \O$, any
limiting point of measures $\tilde\nu^{V_n}_{\w^{(n)}|_{W_{n+1}}}$
is a Gibbs measure. Here $\tilde\nu^{V_n}_{\w^{(n)}|_{W_{n+1}}}$
is a measure on $\O$ such that $\forall n'>n$:
\begin{align*}
&\tilde\nu^{V_n}_{\w^{(n)}|_{W_{n+1}}}\bigg(\{\s\in\O :
\s|_{V_{n'}}=\s_{n'}\}\bigg)\\
=&\left\{ \ba{ll} \nu^{V_n}_{\w^{(n)}|_{W_{n+1}}}(\s_{n'}|_{V_n}),
\ \ \mbox{if} \ \ \s_{n'}|_{V_{n'}\setminus
V_n}=\w^{(n)}|_{V_{n'}\setminus
V_n}\\[2mm] 0, \ \  \ \ \mbox{otherwise}. \ea \right.\,.
\end{align*}

We now recall some basic facts from the theory of von Neumann
algebras. Let $B(H)$ be the algebra of all bounded linear
operators on the Hilbert space $H$ ( over the field of complex
numbers $\C$). A weak (operator) closed $*$-subalgebra ${\mathcal N}$
in $B(H)$ is called {\it von Neumann algebra} if it contains the
identity operator $\e$. By $Proj({\mathcal N})$ we denote the set of
all projections in ${\mathcal N}$. A von Neumann algebra is a {\it
factor} if its center
$$
Z({\mathcal N}:=\{x\in {\mathcal N} : xy=yx, \ \forall y\in {\mathcal N}\}
$$
is trivial, i.e., $Z({\mathcal N})=\{ \l{\e} : \l\in \C\}.$ The von
Neumann algebras are direct sum of the classes I (I$_{n},
n<\infty,$\ I$_{\infty}$), II (II$_{1}$, II$_{\infty})$ and III.
Further, a factor is of only one type among these listed above,
see e.g. \cite{SZ}. An element $x\in{\mathcal N}$ is called positive
if there is an element $y\in {\mathcal N}$ such that $x=y^*y$. A
linear functional $\w$ on ${\mathcal N}$ is called a {\it state} if
$\w(x^*x)\geq 0$ for all $x\in {\mathcal N}$ and  $\w (\e)=1$. A state
$\w$ is said to be normal if
$\w(\sup\limits_{\a}x_{\a})=\sup\limits_{\a}\w(x_{\a})$ for any
bounded increasing net $\{x_{\a}\}$ of positive elements of ${\mathcal
N}$. A state  $\w$ is called {\it trace} (resp. {\it faithful}) if
the condition $\w(xy)=\w(yx)$ holds for all $x,y\in {\mathcal N}$
(resp. if the equality $\w(x^*x)=0$ implies $x=0$).

Let ${\mathcal N}$ be a factor, $\w$ be a faithful normal state on
${\mathcal N}$ and $\s^{\w}_{t}$ be the modular group associated with
$\w$ (see Definition 2.5.15 in \cite{BR1}). We let
$\Gamma(\s^{\w})$ denote the Connes spectrum of the modular group
$\s^{\w}_{t}$ (see Definition 2.2.1 in \cite{BR1}).

For the type III factors, there is a finer classification.

{\bf Definition} (\cite{C}). The type III factor ${\mathcal N}$ is of
type\\
(i) III$_{1}$,  if $\Gamma(\s^{\w}) = \R$;\\
(ii) III$_{\l}$, if $\Gamma(\s^{\w}) =\{n\log \l, n\in {\Z}\},
\quad {\l}\in (0,1)$;\\
(iii) III$_{0}$, if $\Gamma(\s^{\w}) =\{0\}$;\\
see, e.g. \cite{BR1}, \cite{St} for details of von Neumann
algebras and
the modular theory of operator algebras.)

\section{Construction  of Gibbs states for  the
$\l$-model}\label{sec3}

In this section we give a  construction of a special class of
limiting Gibbs measures for the $\l$-model  on the Cayley tree.

Let $h:x\to h_x=(h_{1,x},h_{2,x},...,h_{q-1,x})\in{\R}^{q-1}$ be a
real vector-valued function of $x\in V$. Given $n=1,2,...$
consider the probability measure $\m^{(n)}$ on $\Phi^{V_n}$
defined by
\begin{equation} \mu^{(n)}(\s_n)=Z^{-1}_{n}\exp\{-\b
H(\s_n)+\sum_{x\in W_n}h_x\s(x)\}\,, \label{3.1}
\end{equation}

Here, as before, $\s_n:x\in V_n\to\s_n(x)$ and $Z_n$ is the
corresponding partition function:
\begin{equation*}
Z_n=\sum_{\tilde\s_n\in\Omega_{V_n}}\exp\{\b
H(\tilde\s_n)+\sum_{x\in W_n}h_x\tilde\s(x)\}. \end{equation*}

The consistency conditions for $\m^{(n)}(\s_n), n\geq 1$ are
\begin{equation}
\sum_{\s^{(n)}}\m^{(n)}(\s_{n-1},\s^{(n)})=\m^{(n-1)}(\s_{n-1})\,,
\label{3.2} \end{equation} where $\s^{(n)}=\{\s(x), x\in W_n\}$.

Let $V_1\subset V_2\subset...$ $\cup_{n=1}^{\infty}V_n=V$ and
$\m_1,\m_2,...$ be a sequence of probability measures on
$\Phi^{V_1},\Phi^{V_2},...$ satisfying the consistency condition
(\ref{3.2}). Then, according to the Kolmogorov theorem (see, e.g.
\cite{Sh}), there exists a unique limit Gibbs measure $\m_h$ on
$\O$ such that for every $n=1,2,...$ and $\s_n\in\Phi^{V_n}$ the
equality holds
\begin{equation}
\m\bigg(\{\s|_{V_n}=\s_n\}\bigg)=\m^{(n)}(\s_n)\,. \label{3.3}
\end{equation}

Further we set the basis in ${\R}^{q-1}$ to be
$\z_1,\z_2,...,\z_{q-1}$.

The following statement describes conditions on $h_x$ guaranteeing
the consistency condition of measures $\m^{(n)}(\s_n)$.
\begin{theorem}\label{31}
The measures $\m^{(n)}(\s_n), \ n=1,2,...$ satisfy the consistency
condition  (\ref{3.2}) if and only if for any $x\in V$ the
following equation holds:
\begin{equation}
h'_x=\sum_{y\in S(x)} F(h'_y,\l)\,, \label{3.3*}
\end{equation}

Here, and below $h'_x$ stands for the vector $\frac{q}{q-1}h_x$
and $F:{\R}^{q-1}\to {\R}^{q-1}$ function  is
$F(h;\l)=(F_1(h;\l),...,F_{q-1}(h;\l))$, with
\begin{align}
&F_i(h_1,h_2,...,h_{q-1};\l)\\
=&\log\frac{\sum_{j=1}^{q-1}\exp\{-\b\l(\z_i,\z_j;J)\}\exp
h_j+\exp\{-\b\l(\z_i,\z_q;J)\}}
{\sum_{j=1}^{q-1}\exp\{-\b\l(\z_q,\z_j;J)\}\exp
h_j+\exp\{-\b\l(\s_q,\z_q;J)\}}\,,\nonumber \label{3.5}
\end {align}
$i=1,2,...,q-1$, $h=(h_1,...,h_{q-1}).$
\end{theorem}

The proof uses the same argument as in \cite{MR1},\cite{MR2}.
Denote
\begin{equation*}
{\mathcal D}=\{h=(h_x\in{\R}^{q-1} : x\in  V) :  \  h_x=\sum_{y\in
S(x)} F(h_y,\l),\  \forall  x\in V\}\,.
\end{equation*}

According to  Theorem \ref{31} for  any $h=(h_x, x\in V)\in {\mathcal
D}$ there exists a unique Gibbs measure $\m_h$ which satisfies the
equality (\ref{3.3}).

If the vector-valued function $h^0=(h_x=(0,...,0), x\in V)$ is a
solution, i.e. $h^0\in {\mathcal D}$ then the corresponding Gibbs
measure $\m^{(\l)}_0$ is called {\it the unordered phase} of the
$\l$-model. Since we deal with this unordered phase, we have to
make an assumption which guarantees us the existence of the
unordered phase.

{\bf Assumption A.} {\it For the  considered model the
vector-valued function $h^0=(h_x=(0,0,...,0)$, $x\in V)$ belongs
to ${\mathcal D}$.}

This means that the equation  (\ref{3.3*}) has a solution
$h_x=h_0=0, x\in V$.

According to  Theorem \ref{21} any transformation $S$ of the
group $G_{k+1}$  induces an automorphism $\hat S$  on $V$. By
${\mathcal G}_{k+1}$ we  denote the left group of shifts of $G_{k+1}$.
Any $T\in {\K}_{k+1}$ induces a shift automorphism $\tilde  T:
\O\to\O$  by
\begin{equation*}
(\tilde T\s)(h)=\s(Th), \ \  h\in G_{k+1},\ \s\in \O\,.
\end{equation*}

It is easy  to see that  $\m_0^{(\l)}\circ\tilde T=\m_0^{(\l)}$
for every $\tilde T\in{\K}_{k+1}$.  As mentioned above, the
measure $\m_0^{(\l)}$ has a Markov property (see, \cite{Sp}).

{\bf Assumption B.} {\it We suppose that the measure $\m_0^{(\l)}$
enables a mixing property,  i.e. for any $A,B\in {\mathcal F}$ the
following holds}
\begin{equation}
\lim_{|g|\to\infty}\m_0^{(\l)}(\tilde T_g(A)\cap
B)=\m^{(\l)}_0(A)\m^{(\l)}_0(B)\,. \label{3.8}
\end{equation}

Note that the last condition is satisfied, for example, if  the
phase transition does not occur for the model under 
consideration.\\[2mm]

\section{Diagonal states and  corresponding von Neumann algebras}

Consider  $C^*$-algebra $A=\otimes_{{\G}^k}M_{q} (\C)$, where
$M_{q}({\C})$ is the algebra of $q\times q$ matrices over the
field $\C$ of complex numbers. By  $e_{ij}, \ i,j\in\{1,2,...,q\}$
we denote the basis matrices of the algebra $M_{q}(\C)$. We let
${\rm C}M_{q}(\C)$ denote the commutative subalgebra of
$M_{q}(\C)$ generated by the elements $e_{ii}$ $i=\{1,2,...,q\}$.
We set ${\rm C}A=\otimes_{{\G}^k}{\rm C}M_{q}(C)$. Elements of
commutative algebra ${\rm C}A$ are functions on the space
$\Omega=\{e_{11},...,e_{qq}\}^{{\G}^k}$. Fix a measure $\mu$ on
the measurable space $(\Omega, B)$, where $B$ is the $\s$-algebra
generated by cylindrical subsets of $\Omega$. We construct a state
$\w_{\m}$ on $A$ as follows. Let $P:A\to CA$ be the  conditional
expectation, then the state $\w_{\m}$ is be defined by
$\w_{\m}(x)=\m(P(x))$, $x\in A$, here $\m(P(x))$ means an integral
of a function $P(x)$ under measure $\m$, i.e.
$\m(P(x))=\int_{\O}P(x)(s)d\m(s)$ (see \cite{SV}). The state is
called {\it diagonal}.

By $\w^{(\l)}_0$ we denote the diagonal state generated by the
unordered phase $\mu_0^{(\l)}$. The Markov property implies that
the state   $\w^{(\l)}_0$ is a quantum Markov state (see
\cite{AF}). On a finite dimensional $C^*$-subalgebra
$A_{V_n}=\otimes_{V_n}M_{q}({\C})\subset A$ we rewrite the state
$\w^{(\l)}_0$ as follows
\begin{equation} \w^{(\l)}_0(x)=\frac{tr(e^{\tilde
H(V_n)}x)}{tr(e^{\tilde H(V_n)})}, \ \ x\in A_{V_n}\,, \label{4.1}
\end{equation}
where  $tr$ is the canonical trace on $A_{V_n}$. The term
$\l(\sigma(x)\sigma(y);J)$ in (\ref{3.2}) is given by a diagonal
element of  $M_q(\C)\otimes M_q(\C)$ in the standard basis as
follows
\begin{equation}
\b\l(\sigma(x),\sigma(y);J)= \left ( \ba{cccccc}
B^{(1)}&0&\cdot&\cdot&\cdot&0 \\ 0&B^{(2)}&0&\cdot&\cdot&0\\
\cdot&\cdot&\cdot&\cdot&\cdot&\cdot\\
\cdot&\cdot&\cdot&\cdot&\cdot&\cdot\\
\cdot&\cdot&\cdot&\cdot&\cdot&\cdot\\
0&0&\cdot&\cdot&\cdot&B^{(q)}\\ \ea \right )\,. \label{4.2}
\end{equation}

Here, $B^{(k)}=(b_{ij,k})_{i,j=1}^{q}$, $k=1,dots,q$ are $q\times
q$ matrices, and
\begin{equation}
b_{ij,k}= \left\{ \ba{ll} -\beta\l(\z_k,\z_i;J),   \ i=j,\
i=1,\dots,q\\ 0, \  \  i\neq j \ea \right.\,. \label{4.3}
\end{equation}

Consequently, using (\ref{4.1}) and from (\ref{4.2}),(\ref{4.3})
(cp. \cite{S}, Ch.1, \S 1) the form of Hamiltonian $\tilde H(V_n)$
in the standard basis of $A_{V_n}$ (i.e. under the basis matrices)
is regarded as
\begin{equation*}
\tilde H(V_n)=\sum_{<x,y>\in L_n}\Phi_{x,y}, \ \ \Phi_{x,y}= \left
( \ba{cccccc}
B^{(1)}&0&\cdot&\cdot&\cdot&0 \\
0&B^{(2)}&0&\cdot&\cdot&0\\
\cdot&\cdot&\cdot&\cdot&\cdot&\cdot\\
\cdot&\cdot&\cdot&\cdot&\cdot&\cdot\\
\cdot&\cdot&\cdot&\cdot&\cdot&\cdot\\
0&0&\cdot&\cdot&\cdot&B^{(q)}\\
\ea \right )\,. \end{equation*}

Denote  ${\mathcal M}=\pi_{\w^{(\l)}_0}(A)''$, where
$\pi_{\w^{(\l)}_0}$ is the GNS - representation associated with
the state $\w^{(\l)}_0$ (see Definition 2.3.18 in \cite{BR1}). Our
goal in this section is to determine a type of ${\mathcal M}$.

{\bf Remark.} In \cite{SV} general properties of a representation
associated with diagonal state were studied, but there  concrete
constructions of states were not considered. In \cite{MY} a deep
classification of types of factors generated by quasi-free states
has been obtained. For translation-invariant Markov states the
corresponding type analysis has been made in \cite{GM}. The
investigation of the type of factor arising from translation
invariant or periodic quantum Markov states on the one dimensional
chains is contained in \cite{FM}.

Now we define translations of the $C^*$-algebra $A$. Every
$T\in{\mathcal G}_{k+1}$ induces a translation automorphism $\t_T:A\to
A$ defined by
\begin{equation*}
\t_T(\prod_{x\in V_n}^{\otimes}a_x)= \prod_{x\in
V_n}^{\otimes}a_{T(x)}\,. \end{equation*}

Since measure $\m^{(\l)}_0$ satisfies a mixing property (see
(\ref{3.8})), then we can easily obtain that $\w^{(\l)}_0$ also
satisfies the mixing property under the translations
$\{\t_T\}_{T\in{\K}_{k+1}}$, i.e. for all $a,b\in A$ the equality
holds
\begin{equation*}
\lim_{|g|\to\infty}\w_0^{(\l)}(\t_g(a)b)=\w^{(\l)}_0(a)\w^{(\l)}_0(b)\,.
\end{equation*}
According to Theorem 2.6.10 in \cite{BR1}, the algebra ${\mathcal M}$
is a factor. We note that the modular group of ${\mathcal M}$
associated with $\w^{(\l)}_0$  is defined by
\begin{equation}
\s^{\w^{(\l)}_0}_{t}(x)=\lim_{\L\uparrow\Gamma^k}\exp\{it\tilde
H(\L)\}x\exp\{-it\tilde H(\L)\}, \ \ x\in {\mathcal M}\,, \label{4.4}
\end{equation}
where $\tilde H(\L)=\sum\limits_{<x,y>\in\L}\Phi_{x,y}$. It is
well know that the last limit exists if a suitable norm of the
potential $\tilde H$ is finite (see Theorem 6.2.4.\cite{BR2}).
First of all, we recall the definition of the norm of a potential
$\Psi=\sum_{X\subset \G^k}\Psi(X)$.
\begin{equation*}
\|\Psi\|_{d}=\sum_{n\geq 0}e^{dn}(\sup_{x\in {\G}^k}\sum_{x\in X,
|X|=n+1}\|\Psi(X)\|)\,, \end{equation*} where $d>0$. Here
$\Psi(X)\in A_X=\otimes_XM_q(\mbox{\bf C})$.

Now we compute $\|\tilde H\|_{d}$:
\begin{align*}
\|\tilde H\|_{d}=e^{2d}(\sup_{x\in {\G}^k}\sum_{x\in X,
X=\{u,v\}}\|\Phi_{u,v}\|) =&ke^{2d}\sup_{\{u,v\}\in
L}\|\Phi_{u,v}\|\\
=&ke^{2d}\max_{i,j,k}|\log p_{ij,k}|<\infty\,.
\end{align*}

Hence the norm of $\tilde H$ is finite, therefore the limit in
(\ref{4.4}) exists.

By ${\mathcal M}^{\s}$ one denotes the centralizer of $\w^{(\l)}_0$,
which is defined as
\begin{equation*}
{\mathcal M}^{\s}=\{x\in {\mathcal M} : \s^{\w^{(\l)}_0}_t(x)=x, \ \ \ \
t\in {\R} \}\,. \end{equation*} Since $\w^{(\l)}_0$ is Gibbs
state, according to Proposition 5.3.28 \cite{BR2}, the centralizer
${\mathcal M}^{\s}$  coincides with the set
\begin{equation}
{\mathcal M}_{\w^{(\l)}_0}=\{x\in {\mathcal M} :
\w^{(\l)}_0(xy)=\w^{(\l)}_0(yx), \ \ \ \ y\in {\mathcal M}\}\,,
\label{4.5}
\end{equation}
where we denote by $\w^{(\l)}_0$ also the normal extension of the
state under consideration to all of ${\mathcal M}$.

By $\Pi[n]$ we denote the group of all permutations $\gamma$ of
the set $V_n$ such that
\begin{equation*}
\gamma(x)=x, \ \ \ x\in W_n\,. \end{equation*} Every $\gamma\in
\Pi[n]$ defines an automorphism $\a_{\gamma}: {\mathcal M}\to{\mathcal M}$
by
\begin{equation}
\left. \ba{ll} \a_{\gamma}(\prod_{x\in V_n}^{\otimes}a_x)=
\prod_{x\in V_n}^{\otimes}a_{\gamma(x)}\\[2mm]
\a_{\gamma}\mid_{\bigotimes_{x\notin V_n}M_q({\C})}=id, \ea
\right\}\,, \label{4.6}
\end{equation}
where $id$ is the identity mapping.

Denote
\begin{equation*}
{\mathcal S}_0=\bigcup\{\a_{\gamma} | \gamma\in\Pi[n] \}.\,
\end{equation*}

Simply repeating  the proof of a proposition in \cite{K}, we can
prove the following
\begin{lemma}\label{41}
The group
\begin{equation*}
G_0=\{\a\in {\mathcal S}_0 \mid \  \w^{(\l)}_0(\a(x))=\w^{(\l)}_0(x), \ \ 
x\in {\mathcal M}\}\,, \end{equation*} acts ergodically on ${\mathcal M}$,
i.e. the equality $\a(x)=x$, $\a\in G_0$ implies $x=\theta\e$,
$\theta\in \C$.
\end{lemma}

\begin{lemma}\label{42}
The centralizer  ${\mathcal M}^{\s}$ is a factor of type II$_1$.
\end{lemma}
\begin{pf} From the definition of the automorphism $\a_{\gamma}$
(see (\ref{4.6})), it is easy to see that every automorphism
$\a\in G_0$ is inner, i.e. there exists a unitary $u_{\a}\in {\mathcal
M}$ such that $\a(x)=u_{\a}xu^*_{\a}, x\in {\mathcal M}$. From the
condition $\w^{(\l)}_0\circ\a=\w^{(\l)}_0$ we find
\begin{equation*}
\w^{(\l)}_0(u_{\a}xu^*_{\a})=\w^{(\l)}_0(x), \ \ \ x\in{\mathcal M}\,.
\end{equation*} It follows from (\ref{4.5}) that $u_{\a}\in {\mathcal
M}_{\w^{(\l)}_0}$. According to Lemma \ref{41} the group $G_0$
acts ergodically, this means that the equality $u_{\a}x=xu_{\a}$
for every $\a\in G_0$ implies $x=\theta\e, \theta\in\C$. Hence, we
obtain $\{u_{\a} | \a\in G_0 \}'=\C\e$. Since ${\mathcal
M}^{\s'}\subset\{u_{\a}\}'$ we then get
\begin{equation*}
{\mathcal M}^{\s'}\cap{\mathcal M}=\C\e\,. \end{equation*} In particular
${\mathcal M}^{\s'}\cap{\mathcal M}^{\s}=\C\e$. This means that ${\mathcal
M}^{\s}$ is a factor.
\end{pf}

Now we are able to prove main result of the paper (compare with
the analogous result in \cite{FM}).
\begin{theorem}${}$\\
\label{43}
\begin{itemize}
\item[(i)] If the fraction
\begin{equation*}
\dsf{\l(\z_i,\z_j;J)-\l(\z_m,\z_l;J)}{\l(\z_k,\z_p;J)-\l(\z_u,\z_v;J)}\,,
\end{equation*}
is rational for every $i,j,m,l$, $k,p,u,v\in{1,2,...,q}$, whenever
the denominator is different from 0, then the von Neumann algebra
${\mathcal M}$ associated with the quantum Markov state corresponding
to the unordered phase of $\l$-model (\ref{2.2}) on a Cayley tree
is a factor of type III$_{\theta}$. \item[(ii)] If ${\mathcal M}$ is a
type of III$_1$ factor, then all the fraction cannot be rational.
\end{itemize}
\end{theorem}
\begin{pf} It is known (see Proposition 2.2.2 in \cite{C})
that Connes' spectrum $\Gamma(\a)$  of group of automorphisms
$\a=\{\a_g\}_{g\in G}$ of von Neumann algebra  $M$ has the
following form
\begin{equation}
\Gamma(\a)=\cap\{Sp(\a^e) | e\in Proj(Z(M^{\a})), \ e\neq 0\}\,,
\label{4.7}
\end{equation}
where $\a^e(x)=\a(exe), x\in eMe$ and $Z(M^{\a})$ is the center of
subalgebra
\begin{equation*}
M^{\a}=\{x\in M : \a_g(x)=x, \  \ \ g\in G \}\,. \end{equation*}

Here, $Sp(\a)$ be the Arveson's spectrum of group of automorphisms
$\a$ (see for more details \cite{C},\cite{St}).

By virtue of Lemma \ref{42} we have $Z({\mathcal M}^{\s})=\C\e$. The
equality (\ref{4.7}) implies
$\Gamma(\s^{\w^{(\l)}_0})=Sp(\s^{\w^{(\l)}_0})$.

We now consider the operator $\tilde H(V_n)=\sum\limits_{<x,y>\in
L_n}\Phi_{x,y}$. We let $Sp(\tilde H(V_n))$ denote the spectrum of
the operator $\tilde H(V_n)$. Setting \begin{equation*}
\s^{\w^{(\l)}_0,n}_{t}(x)=\exp\{it\tilde H(V_n)\}x\exp\{-it\tilde
H(V_n)\}, \ \ x\in {\mathcal M}\,, \end{equation*} we obtain
\begin{equation}
Sp(\s^{\w^{(\l)}_0,n})=Sp(\tilde H(V_n))-Sp(\tilde H(V_n))=
\{\l-\m: \l,\m\in Sp(\tilde H(V_n))\}\,. \label{4.8}
\end{equation}

It is clear that  $\b\l(\z_i,\z_k;J)\in Sp(\tilde H(V_n))$,
$\forall i,k\in 1,\dots,q.$ Formula  (\ref{4.8}) then implies that
$Sp(\s^{\w^{(\l)}_0,n})$ is generated by elements of the form
\begin{equation*}
\l(\z_i,\z_j;J)-\l(\z_k,\z_l;J), \ \ i,j,k,l\in 1,\dots q\,.
\end{equation*}
Since
$\frac{\l(\z_i,\z_j;J)-\l(\z_m,\z_l;J)}{\l(\z_k,\z_p;J)-\l(\z_u,\z_v;J)}$
is a rational number, then there is a number $\gamma\in (0,1)$ and
integers $m_{i,j,k,l}\in Z$, ($i,j,k,l\in \{1,2,...,q\}$) such
that
\begin{equation}
\l(\z_i,\z_j;J)-\l(\z_k,\z_l;J)=m_{i,j,k,l}\log\gamma\,.\label{4.9}
\end{equation}

Hence we find that an increasing sequence $\{E(n)\}$ of subsets
$\Z$ such that $E(-n)=-E(n)$ and  $Sp(H(V_n))=\{m\log
\gamma\}_{m\in E(n)}$ is valid. It follows that \begin{equation*}
Sp(\s^{\w^{(\l)}_0})\subset \{n\log \gamma\}_{n\in \Z}\,.
\end{equation*} Hence there exists a positive integer $d\in\Z$ such
that we have
\begin{equation*}
\G(\s^{\w^{(\l)}_0})=\{n\log\gamma^d\}_{n\in \Z}\,.
\end{equation*}
This means that ${\mathcal M}$ is a  factor of type III$_{\theta}$,
$\theta=\gamma^d$. \end{pf}

\section{Applications and examples}

\subsection{Potts model}

We consider the Potts model on the Cayley tree $\G^k$ whose
Hamiltonian is regarded as \begin{equation*}
H(\sigma)=-\sum_{<x,y>}J\delta_{\sigma(x)\sigma(y)}\,,
\end{equation*} where $J\in{\R}$ is a coupling constant, as usual
$<x,y>$ stands for the nearest neighbor vertices and as before
$\s(x)\in \Phi=\{\z_1,\z_2,...,\z_q\}$. Here $\delta $ is the
Kronecker symbol.

Equality (\ref{2.1}) implies that \begin{equation*}
\delta_{\sigma(x)\sigma(y)}=\frac{q-1}{q}
\bigg(\sigma(x)\sigma(y)+\frac{1}{q-1}\bigg)\,, \end{equation*}
for all $ x,y \in V $. The Hamiltonian $H(\sigma)$ is therefore
\begin{equation}
H(\sigma)=-\sum_{<x,y>\in L}J'\sigma(x)\sigma(y)\,, \label{5.1}
\end{equation}
where $ J'=\dsf{q-1}{q}J.$

Hence the $\l$-model is a generalization of the Potts model, that
is in this case the function $\l:\Phi\times\Phi\times {\R}\to
{\R}$ is defined by $\l(x,y;J')=-J'(x,y)$. Here, $x,y\in
{\R}^{q-1}$ and $(x,y)$ stands for the scalar product in
${\R}^{q-1}$. From (\ref{2.1}) it is easy to see that
\begin{equation}
\l(\z_i,\z_j;J)= \left\{ \ba{ll} J, \ \ \textrm{if}\ \ i=j\\[2mm]
-\frac{J}{q-1}, \  \ \textrm{if} \ \  i \neq j \ea\right.\,.
\label{5.2}
\end{equation}

 From (\ref{5.2}) and (\ref{3.3*}) we can check that the assumption
A (see section \ref{sec3}) is valid for the Potts model. So there
exists the unordered phase. From (\ref{5.2}) we can find that the
fraction
$\dsf{\l(\z_i,\z_j;J')-\l(\z_m,\z_l;J')}{\l(\z_k,\z_p;J')-\l(\z_u,\z_v;J')}$,
takes values $\pm 1$ and $0$. So by Theorem \ref{43} a von
Neumann algebra ${\mathcal M}$ is a III$_{\phi}$-factor. From
(\ref{4.9}) we may obtain that
$\phi=\exp\bigg\{\dsf{-J'q}{(q-1)}\bigg\}$. Hence we obtain the
following
\begin{theorem}\label{51}
The von Neumann algebra ${\mathcal M}$ corresponding to the unordered
phase of the Potts model (\ref{5.1}) on a Cayley tree is a factor
of type III$_{\phi^k}$, for some $k\in{\Z}, k>0$, where
$\phi=\exp\bigg\{\dsf{-J'q}{(q-1)}\bigg\}$.
\end{theorem}
{\bf Remark.} If $q=2$ the considered Potts model reduces to the
Ising model, for this model analogous  results were obtained in
\cite{M}. For a class of inhomogeneous Potts model similar result
has been also obtained in \cite{MR1}.

\subsection{Markov random fields}

In this subsection we consider a case when $\l(x,y)$ function is
not symmetric and the corresponding Gibbs measure is a Markov
random field (see \cite{Sp}).

Let $P=(p_{ij})_{i,j=1}^d$ be a stochastic matrix such that
$p_{ij}>0$ for all $i,j\in{1,d}$. Define a function $\l(x,y)$ as
follows:
\begin{equation}
\l(\z_i,\z_j)=-\log p_{ij}\,,\label{6.1}
\end{equation}
for all $i,j\in\{1,\cdots,d\}$. From now on, we will consider the
case $\beta=1$ and $q=d$. It is easy to verify that Assumptions A
and B, for the defined function $\l$, are satisfied. By $\m$ we
denote the corresponding unordered phase of the $\l$-model.
Observe that if the order of the Cayley tree is $k=1$ then the
measure $\m$ is a Markov measure, associated with the stochastic
matrix $P$ (see \cite{Sp}).

By $\w_{\m}$ one  denotes the diagonal state corresponding to the
measure $\m$ on $C^*$-algebra $A=\otimes_{\G^k}M_d(\C)$.
\begin{theorem}\label{52}
Let $P=(p_{ij})_{i,j=1}^d$ be a stochastic matrix such that
$p_{ij}>0$ for all $i,j=1,\dots d$ and at least one element of
this matrix is different from $1/2$, and $\w_{\m}$ be the
corresponding Markov state. If there exist integers $m_{ij}$
$i,j\in\{1,\dots,d\}$, and some number $\a\in(0,1)$ such that
\begin{equation}
\label{con} \frac{p_{11}}{p_{i,j}}=\a^{m_{ij}},
\end{equation}
then $\pi_{\w_{\m}}(A)''$ is a factor of type III$_{\theta}$ for
some $\theta\in(0,1)$.
\end{theorem}

In order to prove this theorem it should be used the rationality
condition of Theorem \ref{43}. Namely, if the condition
(\ref{con}) is satisfied then using (\ref{6.1}) one can see that
the rationality condition holds, so we get the assertion.\\

{\bf Remark.} If all elements of the stochastic matrix $P$ equal
to $1/2$ then the corresponding Markov state $\w_{\m}$ is a trace
and consequently $\pi_{\w_{\m}}(A)''$ is the unique hyperfinite
factor of type
II$_{1}$.

{\bf Remark.} There is a conjecture (see for example, \cite{T})
that every factor associated with GNS representation of a Gibbs
state of a Hamiltonian system having a non-trivial interaction is
of type III$_1$. Theorems \ref{43} and \ref{52} show that the
conjecture is not true even if the Hamiltonian has nearest
neighbour interactions.

\section*{Acknowledgements}

The final part of this work was done within the scheme of NATO-CNR
Fellowship at the Universita di Roma 'Tor Vergata'. We are
grateful to L. Accardi and E. Presutti for useful comments and
observations.

\end{document}